\documentclass [12pt,a4paper,reqno]{amsart}
\usepackage{amsmath,amsthm,amscd}
\usepackage{amsfonts,amssymb}
\allowdisplaybreaks

\newcommand{\be}{\begin{equation}}
\newcommand{\ef}{\end{equation}}
\chardef\bslash=`\\ 

\newcommand{\G}{\Gamma}
\newcommand{\wt}{\widetilde}
\newcommand{\wh}{\widehat}

 \renewcommand{\sectionmark}[1]{}

\newcommand{\iy}{\infty}

\newcommand{\fc}{\frac}
\newcommand{\Te} {Teichm\"{u}ller}

 \usepackage{amsfonts}
\newcommand{\field}[1]{\mathbb{#1}}
\newcommand{\g}{\gamma}

\newcommand{\D}{\field{D}}
\newcommand{\om}{\omega}
\newcommand{\z}{\zeta}
\newcommand{\ov}{\overline}
\newcommand{\vp}{\varphi}
\newcommand{\hC}{\widehat{\field{C}}}
\newcommand{\C}{\field{C}}

\newcommand{\B}{\mathbf{B}}
\newcommand{\T}{\mathbf{T}}

\newcommand{\Belt} {\operatorname{Belt}}

\newcommand{\kp} {\kappa}
\newcommand{\x} {\mathbf x}
\renewcommand{\a} {\alpha}
\newcommand{\ld}{\lambda}

\newcommand{\Fib}{\operatorname{Fib}}

\newcommand{\Teich}{\operatorname{Teich}}

\begin{document}

\title{A link between covering and coefficient theorems for holomorphic functions}
\author{Samuel L. Krushkal}

\bigskip
\begin{abstract}
Recently the author presented a new approach to solving the coefficient problems for various classes of holomorphic functions $f(z) = \sum\limits_0^\iy c_n z^n$, not necessarily univalent.
This approach is based on lifting the given polynomial coefficient functionals
$J(f) = J(c_{m_1}, \dots, c_{m_s}), \ 2 < c_{m_1} < \dots < c_{m_s} < \iy$,
onto the Bers fiber space over universal Teichm\"{u}ller space and applying the analytic and geometric features of Teichm\"{u}ller spaces, especially the Bers isomorphism theorem for Teichm\"{u}ller spaces of punctured Riemann surfaces.

In this paper, we extend this approach to more general classes of functions. In particular, this provides a strengthening of de Branges' theorem solving the Bieberbach conjecture.
\end{abstract}

\date{\today\hskip4mm ({\tt link.tex})}

\maketitle

\bigskip
{\small {\textbf {2020 Mathematics Subject Classification:} Primary: 30C55, 30C62, 30C75; Secondary: 30F60, 32G15, 46G20}

\medskip

\textbf{Key words and phrases:} Taylor coefficients of holomorphic and univalent functions, quasiconformal maps, Teichm\"{u}ller spaces, covering theorems, extremal problems, functional}

\markboth{Samuel Krushkal}{Link between covering and coefficient theorems for holomorphic functions}
\pagestyle{headings}

\bigskip\noindent
{\bf 1. Introductory remarks}.
The holomorphic functionals on various classes of holomorphic (and especially, on univalent) functions
on a disk depending on the Taylor coefficients of these functions play an
important role in various geometric and physical applications of complex analysis, for example, in view of their connection with string theory and with a holomorphic extension of the Virasoro algebra.
These coefficients reflect the fundamental intrinsic features of holomorphy and of conformality.
Estimating these coefficients still remains an important problem in geometric function theory.

Recently the author presented in \cite{Kr5}-\cite{Kr7} a new approach to solving the coefficient problems for various classes of holomorphic functions, not necessarily univalent, which allowed to solve some
long staying problems.
This approach is based on lifting the given polynomial coefficient functionals
  \be\label{1}
J(f) = J(a_{m_1}, \dots, a_{m_s}), \quad 2 < a_{m_1} < \dots < a_{m_s} < \iy,
\end{equation}
onto the Bers fiber space $\Fib \T$ over universal Teichm\"{u}ller space $\T$ and applying the analytic and geometric features of Teichm\"{u}ller spaces, especially the Bers isomorphism theorem for Teichm\"{u}ller spaces of punctured Riemann surfaces.

In fact, one can consider in the similar way the general holomorphic functionals, taking holomorphic functions $J(\mathbf a)$ of $\mathbf a$ in some bounded domains in $\hC^s$ embracing the value
domains of coefficients $a_{m_j}$.

It is important that the space $\T$ is modeled by a strongly rotationally  invariant domain (that is, it contains with any its point representing a univalent function $f(z)$ also the representatives corresponding to pre and post rotations
 \be\label{2}
f(z) \mapsto f_{\a, \beta}(z) = e^{i \beta} f(e^{i\a} z)
\end{equation}
with independent $\a$ and $\beta$ from $[- \pi, \pi]$). Accordingly, the given functional $J$ must be
rotationally homogeneous and hence satisfy
$|J(f_{\a, \beta})| = |J(f)|$.

The general canonical classes of univalent functions are the classes $S$ of functions
$$
f(z) = z + a_2 z^2 + \dots
$$
on the unit disk $\D = \{|z| < 1\}$ and the class $\Sigma$ formed by univalent $\hC$-holomorphic functions with the hydrodynamical normalization
$$
F(z) = z + b_0 + b_1 z^{-1} + b_2 z^{-2} + \dots
$$
on the disk
$\D^* = \{z \in \hC = \C \cup \{\iy\}: \ |z| > 1\}$.
The invertions $F_f(z) = 1/f(1/z)$ of $f \in S$ are nonvanishing (zero free) on $\D^*$.
Let $S_Q$ and $\Sigma_Q$ denote their (dense) subclasses formed by functions admitting quasiconformal
extension to the whole Riemann sphere $\hC$.

We shall also consider the invertions $F_f(z)$ of functions $f(z)$ univalent is a hyperbolic simply connected domain; then $F_f$ is defined in domain
$$
D^{-1} := \{z \in \hC: \ 1/z \in D\}.
$$

Now consider the univalent functions
$$
f(z) = a_1 z + a_2 z^2 + \dots \quad \text{with} \ \  |a_1| = 1
$$
in the unit disk $\D$ continuous on the closed disk $\ov \D$.
The classical Schwarz lemma yields that the boundary curve $L$ of $f(\D)$ must have at least one common point $z_0 = e^{i \beta_0}$. This value is attained at a unique point $z_1 = e^{i \a_0}$; so, using the rotations (2), one obtains a representative of $f$ in the class $S$ (in the general case, the correspondence between the boundary points $z_0$ and $z_1$ must be understand in the sense of the
Carath\'{e}odory prime ends).

Note that all functions $f_{\a, \beta}(z)$ have the same Schwarzian derivative
$$
S_f(z) = \left(\frac{f^{\prime\prime}(z)}{f^\prime(z)}\right)^\prime
- \frac{1}{2} \left(\frac{f^{\prime\prime}(z)}{f^\prime(z)}\right)^2 \quad (z \in \D),
$$
and the chain rule
$$
S_{f_1 \circ f}(z) = (S_{f_1} \circ f) f^\prime(z)^2 + S_f(z).
$$
yields for the M\"{o}bius (fractional linear) maps $w = \g(z)$ the equalities
$$
S_{f_1 \circ \g}(z) = (S_{f_1} \circ \g) \g^\prime(z)^2,
\quad S_{\g \circ f}(z) = S_f(z).
$$
Hence, each $S_f(z)$ can be regarded as a quadratic differential $\vp = S_f(z) dz^2$ on $\D$.
T solution $w(z)$ of the Schwarzian equation $S_w(z) = \vp(z)$ with a given holomorphic $\vp$ iis defined up to a M\"{o}bius transformation of $\hC$.

\bigskip
In this paper, we strengthen some main results from \cite{Kr5}, \cite{Kr6} and extend the indicated approach to more general classes of functions. It will be also shown how the distortion results are naturally connected with the covering theorems for holomorphic functions.

\bigskip\noindent
{\bf 2. Generalized classes of functions}.
Let $L \subset \C$ be an oriented bounded quasicircle separating the points $0$ and $\iy$. Denote its interior and exterior domains by $D$ and $D^*$ (so $0 \in D, \ \iy \in D^*$).
Then, if $\delta_D(z)$ denotes the Euclidean distance of $z$ from the boundary of $D$ and $\ld_D(z) |dz|$ is its hyperbolic metric of Gaussian curvature $-4$, we have
$$
\fc{1}{4} \le \ld_D(z) \delta_D(z) \le 1.
$$
The right hand inequality in (4) follows from the Schwarz lemma and the left from the Koebe one-quoter theorem.
In particular, $\ld_\D(z) = 1/(1 - |z|^2), /  \ld_{\D^*}(z) = 1/(|z|^2 - 1)$.

Consider the classes $S_{Q, \theta}(D)$ of univalent functions in domain $D$ with expansions
$$
f(z) = a_1 z^{-1} + a_2 z^2 + \dots
$$
near the origin $z = 0$ with
$a_1 = e^{i \theta}, \ - \pi \le \theta \le \pi$
admitting quasiconformal extensions onto the complementary domain $D^*$.
The Beltrami coefficients $\mu_f(z) = \partial_{\ov z} f/\partial_z f$ of these extensions run over the unit ball
$$
\Belt(D^*)_1 = \{\mu \in L_\iy(\C): \ \mu(z)|D = 0, \ \ \|\mu\|_\iy  < 1\}.
$$
The Schwarzian derivatives of these functions belong to the complex Banach space $\B(D)$ of hyperbolically bounded holomorphic functions in $D$ with norm
$$
\|\vp\|_{\B(D)} = \sup_D \ \ld_D(z)^2 |\vp(z)|
$$
and run over a bounded domain in the space $\B(D)$. This domain models the universal Teichm\"{u}ller space $\T = \Teich(D)$ with the base point $D$.

It is located in some ball $\{\|\vp\|_\B < C_1\}$ and contains the ball $\{\|\vp\|_\B < C_2\}$
with absolute finite constants $C_1, \ C_2$ depending on $D$.
In this model, the Teichm\"{u}ller spaces of all hyperbolic Riemann surfaces are
contained in $\T$ as its complex submanifolds.

This space corresponds to the equivalence classes $[\mu]$ of $\mu \in \Belt(D^*)_1$, letting $\mu_1, \ \mu_2$
from this ball be equivalent if the corresponding quasiconformal maps $f^{\mu_1}$ and $f^{\mu_2}$ of the sphere $\hC$ coincide on the boundary $L$ (hence, on $\ov D$), and the quotient map
$$
\phi_\T: \ \mu \mapsto [\mu] = S_{f^\mu}
$$
is holomorphic from $L_\iy(D^*)$ to $\B(D)$. This projective map is a split submersion (has the local
holomorphic sections).

Note also that the space $\B(D)$ is dual to the Bergman space $A_1(D)$, a subspace of $L_1(D)$ formed by integrable holomorphic functions (quadratic differentials $\vp(z) dz^2$) on $\D^*$, since every linear functional $l(\vp)$ on $A_1(D)$ is represented in the form
$$
l(\vp) = \langle \psi, \vp \rangle_{\D^*} = \iint\limits_{\D^*} \
(|z|^2 - 1)^2 \ov{\psi(z)} \vp(z) dx dy
$$
with a uniquely determined $\psi \in \B(\D)$.

In a similar fashion, one can consider the classes $\Sigma_{D^*, \theta}$ of univalent functions
in the outer domain $D^*$ with expansions
$$
F(z) =  e^{- i \theta} z + b_0 + b_1 z^{-1} + b_2 z^{-2} + \dots
$$
and using the corresponding ball
$$
\Belt(D)_1 = \{\mu \in L_\iy(\C): \ \mu(z)|D^* = 0, \ \ \|\mu\|_\iy  < 1\},
$$
model the universal Teichm\"{u}ller space $\T = \Teich(D)$ with the base point $D^*$ a s a bounded domain in the space $\B(D^*)$ filled by the Schwarzians $\vp = S_{F^\mu}$. Note that
$\vp(z) = O(z^{-4})$ near $z = \iy$.

\bigskip\noindent
{\bf 3. Statement of results}. Let $G$ be a domain in a complex Banach space $X = \{\mathbf x\}$
and $\chi$ be a holomorphic map from $G$ into the universal Teichm\"{u}ller space $\T = \Teich(D)$ with the base point $D$ modeled as a bounded subdomain of $\B(D)$. Assume that $\chi(G)$ is a (pathwise
connected) submanifold of finite or infinite dimension in $\T$.

The following theorem reveals the intrinsic connection between the covering and distortion features of holomorphic functions. Its special case with $D = \D$ was applied in \cite{Kr7}.

\bigskip\noindent
{\bf Theorem 1}. {\it Let $w(z)$ be a holomorphic univalent solution of the Schwarz differential equation
$$
S_w(z) = \chi(\x)
$$
on $D$ satisfying $w(0) = 0, \ w^\prime(0) = e^{i \theta}$ with the fixed $\theta \in [-\pi, \pi]$
and $\x \in G$ (hence $w(z) = e^{i \theta} z  + \sum_2^\infty a_n z^n$).
Put
 \be\label{3}
|a_{2,\theta}^0| = \sup \{|a_2|: \ S_w \in \chi(G)\},
\end{equation}
and let $w_0(z) = e^{i \theta} z + a_2^0 z^2 + \dots$ be one of the maximizing functions. Then:

(a) For every indicated function $w(z)$ , the image domain $w(D)$ covers entirely the disk
$D_{1/(2 |a_{2,\theta}^0|)} = \{|w| < 1/(2 |a_{2,\theta}^0|)\}$.

The radius value $1/(2 |a_{2, \theta}^0|)$ is sharp for this collection of
functions and fixed $\theta$, and the circle $\{|w| = 1/(2 |a_{2,\theta}^0|)$ contains points
not belonging to $w(\D)$ if and only if $|a_2| = |a_{2,\theta}^0|$
(i.e., when $w$ is one of the maximizing functions).

(b) The inverted functions
$$
W(\zeta) = 1/w(1/\zeta) = e^{i \theta}\zeta - a_2^0 + b_1 \zeta^{-1} + b_2 \zeta^{-2} + \dots
$$
with $\z \in D^{-1}$ map domain $D^{-1}$ onto a domain whose boundary is entirely
contained in the disk} $\{|W + a_{2,\theta}^0| \le |a_{2,\theta}^0|\}$.

{\it (c) Let $D$ be the unit disk $\D$ and
 \be\label{4}
\chi(\x) = c_0 + c_1 z + c_2 z^2 + \dots.
\end{equation}
Let a point $\x_0 \in G$ correspond to the function $w_0(z)$ maximizing the quantity $|a_{2,\theta}^0|$ on $[-\pi, \pi]$. If its Schwarzian
$S_{w_0} = \chi(\x_0) = c_0^0 + c_1^0 + \dots$ satisfies
$$
c_1^0 \ne 0,
$$
then the value $|c_1^0|$ is maximal on this class, i.e.,
$$
|c_1^0| = \{|c_1(w_{\a,\beta})|: \ S_w \in \chi(G); \ \a, \beta \in [- \pi, \pi]\}.
$$
}

\noindent
{\bf Proof}. The proof of parts {\it (a), (b)} follows the classical lines of Koebe's $1/4$
theorem (cf. \cite{Go}).

{\it (a)} Suppose that the point $w = c$ does not belong to the image of $\D$ under the map $w(z)$ defined above. Then $c \ne 0$, and the function
$$
w_1(z) = c w(z)/(c - w(z)) = z + (a_2 + 1/c) z^2 + \dots
$$
also belongs to this class, and hence by (3), $|a_2 +1/c| \le
|a_2^0|$, which implies
$$
|c| \ge 1/(2 |a_2^0|).
$$
The equality holds only when
$$
|a_2 + 1/c| = |1/c| - |a_2| = |a_2^0| \quad \text{and} \ \ |a_2| =
|a_2^0|.
$$

{\it(b)} If a point $\zeta = c$ does not belong to the image $W(\D^*)$, then the function
$$
W_1(z) = 1/[W(1/z) - c] = z + (c + a_2) z^2 + \dots
$$
is holomorphic and univalent in the disk $\D$, and therefore,  $|c + a_2| \le |a_2^0|$.

\bigskip
The assertion {\it (c)} of Theorem 1 is a special case of the following general theorem on
polynomial functionals.

We precede its formulation by some remarks.
To have the rotational symmetry (2), one must take $L = \mathbb S^1$ and hence, $D = \D,
\ D^* = \D^*$, and deal with the corresponding classes $S_Q$ and $\Sigma_Q$ and their appropriate subclasses.

Let $S_{Q,\theta}$ denote the class of univalent functions $w(z) = e^{i \theta} z + a_2 z^2 + \dots$ in $\D$ with $S_w \in \chi(G)$. Consider the union
$$
S^0 = \bigcup_{\theta \in [-\pi,pi]} S_{Q,\theta}
$$
and take its closure $\ov{\mathcal S_Q}$ in the topology of locally uniform convergence on the disk $\D$.
It is compact in this topology. Under appropriate choose of domain $G$, the collection
$\ov{S^0}$ provides the class $S$.

\bigskip\noindent
{\bf Theorem 2}. {\it If a polynomial coefficient functional $J$ of type (1) is strongly rotationally homogeneous on $S^0$ and does not vanish on the function $w_0(z)$ maximizing the coefficient $a_2$ on $\ov{S^0}$, then the maximum of $|J(w)|$ on $\ov{S^0}$ is attained only on this function $w_0(z)$.
}

\bigskip
In particular, all this holds for the functionals
$$
J(f) = a_n(f), \quad n \ge 3.
$$
on the class $S^0$. Then Theorem 2 implies

\bigskip\noindent
{\bf Theorem 3}. {\it For all functions $f(z) \in \ov{S^0}$ and all $n \ge 3$, we have the sharp estimates
$$
\max_{\ov{S^0}} |a_n| = |a_n(w_0)|,
$$
with equality for each $n$ only for the function $w_0(z)$ maximizing $|a_2|$ on $\ov{S^0}$ .   }

\bigskip
This and the next theorems strengthens the de Branges theorem proving the Bieberbach conjecture \cite{DB}.

\bigskip\noindent
{\bf Theorem 4}. {\it If $\ov{\mathcal S_Q} = S$, then the extremal function for any functional
$J(f)$  of type (1) is the Koebe function
$$
\kp_\theta(z) = \fc{z}{(1 - e^{i \theta} z)^2} = z +
\sum\limits_2^\iy n e^{- i(n-1) \theta} z^n, \quad - \pi < \theta \le \pi.
$$
}

Recall that $\kp_\theta$ maps the unit disk onto the complement of the ray
$$
w = -t e^{-i \theta}, \ \ \fc{1}{4} \le t \le \iy.
$$
A special case of Theorem 4 was considered in \cite{Kr5}.

\bigskip\noindent
{\bf 4. A glimpse to Teichm\"{u}ller spaces}.
First we briefly recall some needed results from Teichm\"{u}ller space theory on spaces involved in order to prove Theorem 2; the details can be found, for example, in \cite{Be}, \cite{EKK}, \cite{GL}, \cite{Le}.

This theory is intrinsically connected with univalent functions with quasiconformal extension, and it
is technically more convenient to deal with functions from $\Sigma_Q$. Note also that, in contrast to the generic univalent functions, quasiconformal maps requires three normalization conditions to have uniqueness, compactness, holomorphic dependence on parameters, etc.

\bigskip
The {\bf universal Teichm\"{u}ller space} $\T = \Teich (\D)$ is the space of quasisymmetric homeomorphisms of the unit circle $\mathbb S^1$ factorized by M\"{o}bius maps;  all Teichm\"{u}ller spaces have their biholomorphic copies in $\T$.

The canonical complex Banach structure on $\T$ is defined by factorization of the ball of the Beltrami coefficients (or complex dilatations)
$$
\Belt(\D)_1 = \{\mu \in L_\iy(\C): \ \mu|\D^* = 0, \ \|\mu\| < 1\},
$$
letting $\mu_1, \mu_2 \in \Belt(\D)_1$ be equivalent if the corresponding  quasiconformal maps $w^{\mu_1}, w^{\mu_2}$ (solutions to the Beltrami equation $\partial_{\ov{z}} w = \mu \partial_z w$
with $\mu = \mu_1, \mu_2$) coincide on the unit circle $\mathbb S^1 = \partial \D^*$ (hence, on $\ov{\D^*}$). Such $\mu$ and the corresponding maps $w^\mu$ are called $\T$-{\it equivalent}.

The following important lemma from \cite{Kr7} allows one to use another normalization of
quasiconformally extendable functions.

\bigskip\noindent
{\bf Lemma 1}. {\it For any Beltrami coefficient $\mu \in \Belt(\D^*)_1$ and any $\theta_0 \in [0, 2 \pi]$, there exists a point $z_0 = e^{i \a}$ located on $\mathbb S^1$ so that
$|e^{i \theta_0} - e^{i \a}| < 1$ and such that for any $\theta$ satisfying
$|e^{i \theta} - e^{i \a}| < 1$ the equation
$\partial_{\ov z} w =  \mu(z) \partial_z w$
has a unique homeomorphic solution $w = w^\mu(z)$, which is holomorphic on the unit disk $\D$
and satisfies
 \be\label{5}
w(0) = 0, \quad w^\prime(0) = e^{i \theta}, \quad w(z_0) = z_0.
\end{equation}
Hence, $w^\mu(z)$ is conformal and does not have a pole in $\D$ \ (so
$w^\mu(z_{*}) = \iy$ at some point $z_{*}$ with $|z_{*}| \ge 1$).  }

\bigskip
In particular, this lemma allows one to define the Teichm\"{u}ller spaces using the quasiconformally extendible  univalent functions $w(z)$ in the unit disk $\D$ normalizing these functions by
 \be\label{6}
w(0) = 0, \quad w^\prime(0) = 1, \quad w(1) = 1
\end{equation}
and with more general normalization
$$
w(0) = 0, \quad w^\prime(0) = e^{i \theta}, \quad w(1) = 1.
$$
All such functions are holomorphic in the disk $\D$.

In view of the crucial role of this lemma for the proof of our theorem and for completeness of exposition, we present here its complete proof.

First we establish the assertion of this lemma for $\theta_0 = 0$ corresponding to $z = 1$
and start with the coefficients $\mu$ vanishing in a broader disk $\D_r = \{|z| < r\}, \ r > 1$
(so $w^\mu$ is conformal on $\D_r \Supset \ov{\D}$), and  assume that $\mu \ne \mathbf 0$ (the origin of $\Belt(\D^*)_1$).

Fix $a \in (1, r)$; then $1/a \in \D$.
The generalized Riemann mapping theorem for the Beltrami equation
$\partial_{\ov z} w = \mu(z) \partial_z w$  on $\hC$
implies for a given $\theta \in [0, 2 \pi]$ a homeomorphic solution $\wh w$ to this equation satisfying
 \be\label{7}
\wh w(1/a) = 1/a, \quad \wh w^\prime(1/a) = e^{i \theta} \ , \quad \wh w(\iy) = \iy.
\end{equation}

Since, by the classical Schwarz lemma, for any holomorphic map
$g: \D \to \D$ and any point $z_0 \in \D$,
$$
|g^\prime(z_0)| \le (1- |g(z_0)|^2)/(1 - |z_0|^2)
$$
with equality only for appropriate M\"{o}bius automorphism of $\D$, the above  normalization (6) and the assumption on $\mu$ yield for the constructed map $\wh w(z)$ that the image $\wh w(\D)$ does not cover $\D$, and thus either $\wh w(\D)$ is a proper subdomain of $\D$ or it also contains  the points $z$ with $|z| > 1$ outer for $\D$.

Applying this to suitable rotated map
$\wh w_\a(z) =  e^{-i \a} \wh w(e^{i \a} z)$
having Beltrami coefficient $\mu_\a(z) = \mu(e^{i \a} z) e^{2 i \a}$, one obtains that
domain $\wh w_\a(\D)$ does not contain simultaneously both distinguished points
$a$ and $1/a$, at least sufficiently close to $1$ (and the same is valid for the points
$a^\prime \in \D^*$ close to $a$).
Now consider the compositions
$$
w_{a,\a}(z) = \gamma_a^{-1} \circ \wh w_\a \circ \gamma_a(z),
$$
with
$$
\g_a(z) = (1 - a z)/(z - a), \quad a \in \D^*,
$$
The Beltrami coefficient of $w_{a,\a}$ equals $\gamma_{a,*} \mu$. Since
$$
\gamma_a(\infty) = - a, \quad \gamma_a(a) = \infty, \quad \gamma_a(0) = - 1/a,
\quad \gamma_a(1/a) = 0
$$
and accordingly,
$$
\g_a^{-1}(\iy) = a, \ \g_a^{-1}(- a) = \iy, \ \g_a^{-1}(0) = 1/a, \ \g_a^{-1}(- 1/a) = 0,
$$
the map $w_{a,\a}$  satisfies
 \be\label{8}
w_{a,\a}(0) = 0, \quad w_{a,\a}^\prime(0) = \wh w_\a^\prime(1/a) = e^{i \theta}, \quad
w_{a,\a} (a) = \gamma_a^{-1} \circ \wh w_\a \circ \gamma_a(a) = a.
\end{equation}

In view of our assumptions on the map $\wh w_\a$, the point $\wh w_\a^{-1}(a)$ does not lie in the unit disk $\D$, which implies that the function $w_{a,\a}$ is holomorphic in this disk.

\bigskip
Now we investigate the limit process as $a \to 1$.
Any from the constructed maps $w_{a,a}$ is represented as a composition of a fixed solution $\wh w$ to the equation
$\partial_{\ov z} w = \mu(z) \partial_z w$ subject to (7) and some  M\"{o}bius maps $\wh \gamma_a$. The first two conditions in (7) imply that the restrictions of these $\wh \gamma_a$ to $\wh w(\D_r)$ form a
(sequentially) compact set of $\wh \gamma_a$ in the topology of convergence in the spherical metric on $\hC$. Letting $a \to 1$, one
obtains in the limit the map $\wh \gamma_1(z) = \lim\limits_{a\to 1} \wh \gamma_a(z)$,
which also is a non-degenerate (nonconstant) M\"{o}bius map. Accordingly,
$$
\lim\limits_{a\to 1} w_{a,\a}(z) = \wh \gamma_1 \circ \wh w_\a \circ \g_1(z) =: \wh w_1(z),
$$
and this map satisfies (8) with $a = 1$, which is equivalent to (9).

Note that the relations (8) do not depend on $r$ and that the normalization (9) also holds for the inverse rotation $e^{i \a} \wh w_1(e^{-i \a} z)$ of the limit function. Letting
$$
w(z) = e^{i \a} \wh w_1(e^{-i \a} z), \quad z_0 = e^{i \a},
$$
one obtains a weakened assertion of Lemma 1 for all Beltrami coefficients $\mu \ne \mathbf 0$ supported in the disk $\D_r^* = \{|z| > r\}$ with $r > 1$ (without a restriction for $|e^{i\theta} - z_0|$).

\bigskip
To extend the obtained result to arbitrary $\mu \in \Belt(\D^*)_1 \ (\mu \ne \mathbf 0$), we pass to the truncated coefficients
$$
\mu_r(z) = \mu(r z), \quad |z| > 1,
$$
which are equal to zero on the disk $\{|z| < r\}$. The compactness properties of the $k$-quasiconformal families (i.e., with $\|\mu\|_\iy \le k < 1$)
imply the convergence of maps $w^{\mu_r}(z)$ normalized by (5) to $w^\mu$(z) as $r \to 1$
in the spherical metric on $\hC$ (and hence everywhere on $\hC$). Accordingly, one must now take
$z_0 = \lim\limits_{r\to 1} w^{\mu_r}(e^{i \a})$.

It remains to estimate the lower bond for $|e^{i\theta} - z_0|$ and consider the case $\mu(z) \equiv 0$ omitted above. We consider for this the homotopy functions
$$
w_t(z) w^{\mu_t}(z) = c(t) w(t z): \ \C \times (0, 1) \to \C
$$
with $\mu_t(z) = \mu(t z), \ 0 < t < 1$. The factor $c(t)$ is determined by normalization (8)
and is a fractional linear function of $t$.

For any $t$, the points $w_t(a)$ and $w_t(1/a)$ as well as $0$ and $\iy$ are separated by the quasicircle $w_t(\mathbb S^1)$.
Thus, arguing similar to above, one obtains for any of these $w_t$ the same point $z_0 = e^{i \a}$,
and this also holds for $t = 0$.

It remains to observe that if $\mu \to \mathbf 0$ in $L_\iy$ norm (or even $\mu(z) \to 0$ almost everywhere in $\D^*$), then the
corresponding limit map $w(z)$ satisfying (5) must be an elliptic fractional linear transformation with fixed points $0$ and $z_0$; hence,
$$
\frac{w - z_0}{w} = e^{- i\theta} \ \frac{z - z_0}{z},
$$
which implies
  \be\label{9}
w = \fc{e^{i \theta} z}{(1 - e^{- i \theta})z_0^{-1} z + 1}.
\end{equation}
The value $w(z_{*}) = \iy$ occurs when
$$
(1 - e^{-i \theta}) z_0^{-1} z_{*} + 1  = 0.
$$
Then $|z_{*}| = 1/|e^{i \theta} - 1|$, and hence, $|z_{*}| < 1$ if $|e^{i \theta} - 1| > 1$,
what is excluded by assumption.

This implies the assertion of Lemma 1 for $\theta_0 = 0$ and Beltrami coefficient $\mu_\a$.
To get it for $\mu$, one must conjugate $w^{\mu_\a}$ by rotation
$z \mapsto e^{- i \a} z$, replacing the fixed point $z_0 = 1$ by $e^{i \a}$.

Similarly, the case of arbitrary $\theta_0$ satisfying $|e^{i \theta_0} - e^{i \a}| < 1$ is
reduced to the previous one by compositions of $w$ with pre and post rotations
$z \mapsto e^{i \theta_0} z$, completing the proof of the lemma.

\bigskip
It follows from Lemma 1 and from its proof that for any fixed $\theta_0 \in [- \pi, \pi]$ there is a point $z_0 = e^{i \a_0} \in \mathbb S^1$ such  that for all $\theta$ with $|e^{i \theta} - z_0| < 1$
any two Beltrami coefficients $\mu_1, \mu_2 \in \Belt(\D^*)_1$ generate quasiconformal maps $w^{\mu_1}$ and $w^{\mu_2}$ normalized by (5) (hence, having the same fixed point $z_0$),
unless these maps are conjugated by a rotation, or equivalently, $\mu_2(z) = \mu_1(e^{i\a} z) e^{-2i\a}$ with some $\a \in[- \pi, \pi]$.

\bigskip
The proof of Theorem 2 also involves other Teichm\"{u}ller spaces. The corresponding space $\T_1 = \Teich(\D_{*})$ {\bf for the punctured disk} $\D_{*} = \D \setminus \{0\}$ is formed by classes $[\mu]_{\T_1}$ of $\T_1$-{\bf equivalent} Beltrami coefficients $\mu \in \Belt(\D)_1$ so that the corresponding quasiconformal automorphisms $w^\mu$ of the unit disk coincide on both boundary components (unit circle $\mathbb S^1$ and the puncture $z = 0$) and are homotopic on $\D \setminus \{0\}$.
This space can be endowed with a canonical complex structure of a complex Banach manifold
and embedded into $\T$ using uniformization of $\D_{*}$ by a cyclic parabolic Fuchsian
group acting discontinuously on $\D$ and $\D^*$. The functions $\mu \in L_\iy(\D)$ are lifted to
$\D$ as the Beltrami measurable $(-1, 1)$-forms  $\wt \mu d\ov{z}/dz$ in $\D$ with respect to
$\G$, i.e., via $(\wt \mu \circ \g) \ov{\g^\prime}/\g^\prime = \wt \mu, \ \g \in \G$,
forming the Banach space $L_\iy(\D, \G)$.

We extend these $\wt \mu$ by zero to $\D^*$ and consider the unit ball $\Belt(\D, \G)_1$ of
$L_\iy(\D, \G)$. Then the corresponding Schwarzians $S_{w^{\wt \mu}|\D^*}$ belong to $\T$.
Moreover, $\T_1$ is canonically isomorphic to the subspace $\T(\G) = \T \cap \B(\G)$,
where $\B(\G)$ consists of elements $\vp \in \B$ satisfying $(\vp \circ \g) (\g^\prime)^2 = \vp$
in $\D^*$ for all $\g \in \G$.

\bigskip
Due to {\bf the Bers isomorphism theorem}, {\it the space $\T_1$ is biholomorphically isomorphic
to the Bers fiber space
$$
\Fib(\T) = \{(\phi_\T(\mu), z) \in \T \times \C: \ \mu \in
\Belt(\D)_1, \ z \in w^\mu(\D)\}
$$
over the universal Teichm\"{u}ller space $\T$ with holomorphic projection $\pi(\psi, z) = \psi$} (see \cite{Be}).

This fiber space is a bounded hyperbolic domain in $\B \times \C$ and represents the collection of domains $D_\mu = w^\mu(\D)$ as a holomorphic family over the space $\T$. For every $z \in \D$,  its
orbit $w^\mu(z)$ in $\T_1$ is a holomorphic curve over $\T$.

The indicated isomorphism between $\T_1$ and $\Fib(\T)$ is induced by the inclusion map \linebreak
$j: \ \D_{*} \hookrightarrow \D$ forgetting the puncture at the origin via
 \be\label{10}
\mu \mapsto (S_{w^{\mu_1}}, w^{\mu_1}(0)) \quad \text{with} \ \
\mu_1 = j_{*} \mu := (\mu \circ j_0) \ov{j_0^\prime}/j_0^\prime,
\end{equation}
where $j_0$ is the lift of $j$ to $\D$.

The Bers theorem is valid for Teichm\"{u}ller spaces $\T(X_0 \setminus \{x_0\})$ of all punctured hyperbolic Riemann surfaces $X_0 \setminus \{x_0\}$; we use only its special case.

\bigskip
The spaces $\T$ and $\T_1$ can be weakly (in the topology generated by the spherical metric on $\hC$) approximate by finite dimensional Teichm\"{u}ller spaces $\T(0, n)$ of punctured spheres (Riemann surfaces of genus zero)
$$
X_{\mathbf z} = \hC \setminus \{0, 1, z_1 \dots, z_{n-3}, \iy\}
$$
defined by ordered $n$-tuples $\mathbf z = (0, 1, z_1, \dots, z_{n-3}, \iy), \ n > 4$ with distinct
$z_j \in \C \setminus \{0, 1\}$ (the details see, e.g., in \cite{Kr5}).

Another canonical model of $\T(0, n) $ is obtained again using the uniformization. This space
is biholomorphic to a bounded domain in the complex Euclidean space $\C^{n-3}$.

\bigskip\noindent
{\bf 5. Proof of Theorem 2}. The proof given below follows the lines of \cite{Kr5} and involves lifting the functional $J$ onto the spaces $\T$ and $\T_1$.

\noindent
{\it(a)} \ Taking into account Lemma 1, we distinguish the subclasses $S_{z_0,\theta}$ consisting of $f \in S_{Q,\theta}$ with fix point at $z_0 \in \mathbb S^1$; their union
$$
\bigcup_{z_0 \in \mathbb S^1, \theta \in [-\pi,pi]} S_{z_0,\theta}
$$
coincides with the class $\wh S^0$ introduced in \textbf{2.2}.
Further, we pass to the inverted functions $F_f(z) = 1/f(1/z)$ for $f \in \wh S^0$, which form the corresponding classes $\Sigma_{z_0,\theta}$ of nonvanishing  univalent functions on the disk $\D^*$ with expansions
$$
F(z) =  e^{- i \theta} z + b_0 + b_1 z^{-1} + b_2 z^{-2} + \dots,
\quad  F(z_0) = z_0,
$$
and
$$
\Sigma^0 = \bigcup_{z_0,\theta} \Sigma_{z_0,\theta}.
$$
The coefficients $a_n$ of $f(z)$ and the corresponding coefficients $b_j$ of $F_f(z)$ are related by
$$
b_0 + e^{2i \theta} a_2 = 0, \quad b_n + \sum \limits_{j=1}^{n}
\epsilon_{n,j}  b_{n-j} a_{j+1} + \epsilon_{n+2,0} a_{n+2} = 0,
\quad n = 1, 2, ... \ ,
$$
where $\epsilon_{n,j}$ are the entire powers of $e^{i \theta}$. This
successively implies the representations of $a_n$ by $b_j$ via
 \be\label{11}
a_n = (- 1)^{n-1} \epsilon_{n-1,0}  b_0^{n-1} - (- 1)^{n-1} (n - 2)
\epsilon_{1,n-3} b_1 b_0^{n-3} + \text{lower terms with respect to} \ b_0.
\end{equation}
This transforms the initial functional (1) into a coefficient functional $\wt J(F^\mu)$
on $\Sigma^0$ depending on the corresponding coefficients $b_j$. This dependence is holomorphic
from the Beltrami coefficients $\mu_W \in \Belt(\D)_1$ and from the Schwarzians $S_{W^\mu}$.

Accordingly, we model the universal Teichm\"{u}ller space $\T$ by a domain in the space $\B = \B(\D^*)$ formed by the Schwarzians $S_F^\mu$. Thereby, both functionals  $\wt J(F)$ and $J(f)$ are lifted holomorphically onto the space $\T$.

To lift $J$ onto the covering space $\T_1$, we again pass to functional $\wh J(\mu) = \wt J(F^\mu)$
lifting $J$ onto the ball $\Belt(\D)_1$ and apply the $\T_1$-equivalence, i.e., the quotient map
$$
\phi_{\T_1}: \ \Belt(\D)_1 \to \T_1, \quad \mu \to [\mu]_{\T_1}.
$$
Thereby the functional $\wt J(F^\mu)$ is pushed down to a bounded holomorphic functional $\mathcal J$
on the space $\T_1$ with the same range domain.

Using the Bers isomorphism theorem, one can regard the points of the space $\T_1$ as the pairs
$X_{F^\mu} = (S_{F^\mu}, F^\mu(0))$, where $\mu \in \Belt(\D)_1$ obey $\T_1$-equivalence (hence, also $\T$-equivalence). Note that since the coefficients $b_0, \ b_1, \dots$ of $F^\mu \in \Sigma_\theta$   are uniquely determined by its Schwarzian $S_{F^\mu}$, the values of $\mathcal J$ in the
points $X_1, \ X_2 \in \T_1$ with $\iota_1(X_1) = \iota_1(X_2)$ are equal.

In result, we get on $\T_1 = \Fib (\T)$ the holomorphic functional
$$
\mathcal J(X_{F^\mu}) = \mathcal J(S_{F^\mu}, \ t), \quad t = F^\mu(0),
$$
and must investigate the restriction of plurisubharmonic functional $|\mathcal J(S_{F^\mu}, t)|$ to
the image of domain $\chi(G)$ in $\Fib (\T)$. .

Note that by the part {\it (b)} of Theorem 1, the boundary of domain $W^\mu(\D^*)$ under any function
$W^\mu(z) \in \Sigma^0$ is located in the disk $\{|W - b_0| \le |a_2^0| \}$  with $|a_2^0|$
determined by (3). For all such $W^\mu$, the variable $t$ in the representation (11) runs over some
subdomain $D_1$ in the disk $\D_4 = \{|t| < 4\}$ containing the origin (this subdomain depends on $z_1$).
Since the functional $J$ is rotationally invariant, this subdomain $D_1$ is a disk $\D_{\a_1}$ of some radius $\a_1 \le 2 |a_2^0|$.

We define on this domain the function
 \be\label{12}
\wt u_1(t) = \sup_{S_{W^\mu}} \mathcal J(S_{W^\mu}, t).
\end{equation}
taking the supremum over all $S_{W^\mu} \in \T$ admissible for a given $t = W^\mu(0) \in D_{\a_n}$, that means over the pairs $(S_{W^\mu}, t) \in \Fib(\T)$ with a fixed $t$ and pass to the upper-semicontinuous
regularization
$$
u_1(t) = \limsup\limits_{t^\prime \to t} \wt u_1(t^\prime).
$$

\noindent
{\it (b)} \ Now the crucial step in the proof of Theorem 2 is to establish that the function (12) inherits subharmonicity. In fact, we have much more.

Namely, selecting on the unit circle a dense subset
$$
\mathbf e = \{z_1, z_2, \dots, z_n, \dots\}, \quad z_1 = e^{i \theta_1},
$$
and repeating successively for the above construction with fix points $z_1, z_2, \dots$, one obtains  similar
to (12) the corresponding functions $u_1(t), u_2(t), \dots$. Let $u(t)$ be their upper envelope
$\sup_n u_n(t)$ followed by its upper semicontinuous reglarization.
Then we have:

\bigskip\noindent
{\bf Lemma 2}. {\it The functions $u(t)$ is logarithmically subharmonic in some disk $\D_|a$ with
$\a \le 2|a_2^0|$.  }

\bigskip
In view of assumptions on the domain $\chi(G)$, this lemma is a consequence of the general basic Lemma 1 from \cite{Kr5}. The indicated Lemma 1 in \cite{Kr5} establishes subharmonicity of the maximal function
$u(t)$ taking the supremum of $\mathcal J(S_{F^\mu},t))$ over all $S_{F^\mu}$ running over the whole space $\T$. Now we need take this supremum only over $S_{F^\mu} \in \chi(G)$. The restriction of $\mathcal J(S_{F^\mu},t))$ to the image of $\chi(G)$ in $\Fib(\T)$ also is plurisubharmonic.
Note that the proof of the indicated Lemma 1 in \cite{Kr5} essentially involves the strong rotational homogeneity (2) of the initial functional $J(f)$.

\bigskip\noindent
{\it (c)} \ Finally, we have to establish the range domain of $F^\mu(0)$ for $S_{F^\mu}$ running over
$\chi(G)$ and describe the boundary points of this domain.

The assumptions on $\chi(G)$ and the features of construction of the Bers fiber space
$\Fib(\T)$ imply that image of $\chi(G)$ in the space $\Fib(\T)$ also is a connected submanifold
covering $\chi(G)$. This allows one to apply to this image the same arguments as in \cite{Kr5}, \cite{Kr6}, applied there to the whole space $\Fib(\T)$.

It follows from the previous step that the indicated domain is rotationally symmetric and connected,
hence a disk $\D_\a = \{|t| < \a\}$ of some radius $\a \le 2 |a_2^0|$. Theorem 1 (parts {\it (a)} and {\it (b)}) yields that the extremal value of $\a$ equals $2|a_2^0|$.
The maximum of $|J|$ must  be attained on the boundary circle $|t| = 2 |a_2^0|$, corresponding by
Theorem 1 the function $w_0$, and the assertion of  Theorem 2 follows.

\bigskip\noindent
{\bf 6. Remarks on rotational symmetry}.
The assumption on strong rotational symmetry (2) (naturally connected with Schwarz's lemma)
is crucial for applying the above approach, though, for example, Koebe's function $\kappa_\theta$
also is extremal for some functionals which do not have such symmetry.
For example, this holds for Zalcman's functional
$$
Z_n(f) = a_{2n-1} - a_n^2,
$$
which admits only a weakened symmetry $|Z_n(f_\a)| = |Z_n(f)|$ if
$f_\a(z) = e^{- i \a} f(e^{i \a} z)$.

The well-known Zalcman conjecture stated in beginning of 1960's and implying the Bieberbach conjecture says that for any $f \in S$ and all $n \ge 3$,
$$
|a_n^2 - a_{2n-1}| \le (n-1)^2,
$$
with equality only for the Koebe function $\kp_\theta(z)$.

Zalcman's conjecture has been proved by the author for $n = 3, 4, 5, 6$  in \cite{Kr2}, \cite{Kr4} (see also \cite{HL}, \cite{TTV}) and remains open for $n \ge 7$.
\footnote{There are also some special results on this conjecture obtained in \cite{EV}, \cite{OT} for $S$, and this conjecture was proven for some subclasses of $S$.}

On the other hand, the extremals of the functional
 \be\label{13}
J(f) = a_2^3 - 2 a_2 (a_2^2 - a_3) - a_4 + M (a_2^2 - a_3)^2,
\end{equation}
(equal to $b_2 + M b_1^2$ in terms of coefficients of the inverted function  $F_f(z) = 1/f(1/z), \ |z| > 1$) with small $|M|$, are different from the Koebe function.

This follows from the well-known Schiffer's result \cite{Sc}, that the coefficient $b_2$ is sharply estimated on $\Sigma$ by $|b_2| \le 2/3$; clearly, its extremal function is different from
$$
F_\theta(z) = 1/\kp_\theta(1/z) = z - 2 e^{i \theta} + e^{2i \theta} z^{-1}.
$$
The same is valid for the Fekete-Szeg\"{o} functional
 \be\label{14}
J(f) = a_3 - \ld \ a_2^2, \quad 0 < \ld < 1.
\end{equation}
The corresponding to (13), (14) functionals $\wt J(F)$ on $\Sigma$ contain a separated free term
$b_0^j = a_2^j$, whose role is twofold: it contributes the distortion of $\wt J(F)$ and
simultaneously determines the normalization of $F_f$ (and of $f$).
This term does not be controlled by elements of the space $A_1(\D)$.

\bigskip\noindent
{\bf 7. Example}. Consider the balls $B_r = \{\|\vp\|_2 < r\}$ in the Hilbert space $A_2(\D)$ of the
square integrable holomorphic functions
$$
\vp(z) = c_0 + c_1 z + c_2 z^2 + \dots, \quad |z| < 1,
$$
with norm
$$
\|\vp\|_2 = \Bigl(\fc{1}{2\pi} \iint\limits_\D |\om|^2 dx dy\Bigr)^{1/2}.
$$
The Taylor coefficients $c_n$ of these functions coincide with their Fourier coefficients with respect
to the orthonormal system
 \be\label{15}
\vp_n(z) = \sqrt{\fc{n + 1}{\pi}} \ z ^ n, \quad n = 0, 1, 2, \dots , \quad \text{on} \ \ \D,
\end{equation}
and
$$
\sum\limits_0^\iy |c_n|^2 = \|\vp\|_2.
$$
Applying the Schwarz inequality, one obtains
$$
\fc{1}{\pi} \iint\limits_\D |\vp| dx dy \le \fc{1}{\sqrt{\pi}}
\Bigl(\iint\limits_\D |\vp|^2 dx dy\Bigr)^{1/2} = \sqrt{2}   \pi \ \|\vp\|_2,
$$
which together with the inequality
$$
\|\vp\|_\B \le \iint_\D |\vp| dx dy
$$
(following from the mean value inequality for integrable holomorphic functions in a domain) implies
that all functions $\vp \in A_2(\D)$ belong to the space $\B$, and their $\B$-norm is estimated by
 \be\label{16}
\|\vp\|_\B \le \sqrt{2}\pi \  \|\vp\|_2.
\end{equation}
Hence, these functions (more precisely, the corresponding quadratic differentials $\vp(z) dz^2$) can be regarded as the Schwarzian derivatives of locally univalent functions in $\D$.

Now, take the domain $G = B_{r_0}$ with
$$
r_0 = \sqrt{2}/\pi.
$$
By (16), all $\vp \in G$ have the $\B$-norm less than $2$, and by the Ahlfors-Weill theorem \cite{AW}
are the Schwarzian derivatives of univalent functions $w(z)$ on the disk $\D$ having canonical quasiconformal extensions to $\hC$.
Hence, every coefficient $c_n$ is represented as a polynomial from the initial coefficients
$a_2, \dots, a_{n+1}$ of $w(z)$.
The maximal value of each $|c_n|$ equals $1$ and is attained on $G$ on the function $\vp_n$.

Theorems 1 and 2 yield that the univalent function $w_0(z)$ with $S_{w_0}= \vp_1$ maximizes simultaneously both coefficients $a_2$ and $c_1$ on this class, and this function also is extremal
for all $a_n, \ n \ge 3$.

To determine explicitly the Taylor coefficients $a_n$ of univalent $w(z)$ whose Schwarzians run over the closed ball $\ov{B_{r_0}}$, one has to find the ratios $w = \eta_1/\eta_2$ of two independent solutions of the linear differential equation
$$
2 \eta^{\prime\prime}(z) + \vp(z) \eta(z) = 0
$$
with a given holomorphic $\vp$ (which is equivalent to solving the Schwarzian equation $S_w(z) = \vp(z)$)
subject to $w(0) = 0, \ w^\prime(0) =1 $.
The arising linear equations are of the form
$$
\eta^{\prime\prime} + c z^p \ \eta = 0
$$
and closely relate to special functions (see, e.g., \cite{Ka}). For the extremal function $w_0$
distinguished above, we have the equation
$$
\eta^{\prime\prime} + r_0 z \ \eta = 0;
$$
its solutions are represented by linear combinations of cylindrical functions.

\bigskip
Note also that, due to the existence theorem from \cite{Kr3}, for every given $n \in \mathbb N$ and
every bounded holomorphic function $\vp(z) = c_0 + c_1 z + \dots \in A_2(\D)$, which is not a polynomial of degree $n_1 \le n$, there exists a quasiconformal deformation of the extended plane $\hC$ conformal on $\vp(\D)$, which perturbs arbitrarily (in restricted limits) the coefficients
$c_0, c_1, \dots, c_{n-1}$ but preserves the $L_2$-norm of $\vp$.
The maximizing functions (15) do not admit such deformations.

\medskip
{\small\em{ \leftline{Department of Mathematics, Bar-Ilan
University, 5290002 Ramat-Gan, Israel} \leftline {Department of Mathematics,
University of Virginia, Charlottesville, VA 22904-4137, USA}}}


\bigskip
\bigskip
\begin{thebibliography}{EKK}
{\small

\bibitem{AW}
L.V. Ahlfors and G. Weill, {\it A uniqueness theorem for Beltrami equations}, Proc. Amer. Math. Soc. \textbf{13} (1962), 975-978.

\bibitem{Be}
L. Bers, {\it Fiber spaces over Teichm\"{u}ller spaces}, Acta Math. \textbf{130} (1973), 89-126.

\bibitem{DB}
L. de Branges, {\it A proof of the Bieberbach conjecture}, Acta Math. \textbf{154}
(1985), 137-152.

\bibitem{EKK}
C.J. Earle, I. Kra and S. L. Krushkal,
{\it Holomorphic motions and \Te \ spaces},
Trans. Amer. Math. Soc.  \textbf{343} (1994), 927-948.

\bibitem{EV}
I.Eframidis and D. Vukoti\'{c}, {\it Application of Livingston-type inequalities to the generalized Zalcman functional}, Math. Nachr. \textbf{291} (2018), 1502-1513.

\bibitem{GL}
F.P. Gardiner and N. Lakic, {\it Quasiconformal \Te \ Theory}, Amer. Math. Soc.,
Providence, RI, 2000.

\bibitem{Go}
G.M. Goluzin, {\it Geometric Theory of Functions of Complex Variables}, Transl. of Math. Monographs,
vol. 26, Amer. Math. Soc., Providence, RI, 1969.

\bibitem{HL}
W.K. Hayman and E.F. Lingham, {\it Research Problems in Function Theory - Fiftieth Aniversary Edition}, Problems Books in Mathematics, Springer, 2019.

\bibitem{Ka}
E. Kamke, {\it Differentialgleichungen: L\"{o}sungsmethoden und L\"{o}sungen I. Gew\"{o}hnliche
Differentialgleichungen}, 6. Aufl., Leipzig, 1959.

\bibitem{Kr1}
S.L. Krushkal, {\it Quasiconformal Mappings and Riemann Surfaces}, Wiley, New York, 1979.

\bibitem{Kr2}
S.L. Krushkal, {\it Univalent functions and holomorphic motions},
J. Analyse Math. \textbf{66} (1995), 253-275.

\bibitem{Kr3}
S.L. Krushkal, {\it Quasiconformal maps decreasing $L_p$ norm}, Siberian Math. J. \textbf{41} (2000), 884-888.

\bibitem{Kr4}
S.L. Krushkal, {\it Proof of the Zalcman conjecture for initial
coefficients}, Georgian Math. J. \textbf{17} (2010), 663-681;
Corrigendum \textbf{19} (2012), 777.

\bibitem{Kr5}
S.L. Krushkal {\it Teichm\"{u}ller spaces and coefficient problems for univalent holomorphic functions}, Analysis and Mathematical Physics \textbf{10} (2020), no. 4.
https://doi.org/10.1007/s13324-020-00395-y

\bibitem{Kr6}
S.L. Krushkal, {\it Teichm\"{u}ller space theory and classical problems of geometric function theory},
J. Math. Sci. \textbf{258} (3) (2021), 276-289.
DOI: 10.1007/s10958-021-05546-5

\bibitem{Kr7}
S.L. Krushkal, {\it Two coefficient conjectures for nonvanishing Hardy functions, I}, J. Math. Sci. \textbf{268} (2) (2022), 199-221.
DOI: 10.1007/s10958-022-06192-1

\bibitem{Le}
O. Lehto, {\it Univalent Functions and Teichm\"uller Spaces},
Springer-Verlag, New York, 1987.

\bibitem{OT}
M. Obradovi\'{c} and N. Tuneski, {\it Certain properties of the class of
univalent functions with real coefficients}, arXiv: 2112.15449v1 [mathCV], 2021.

\bibitem{Sc}
M. Schiffer, {\it Sur un probleme d'extremum de la representation conforme}, Bull. Soc. Math.
France \textbf{66} (1938), 48-56.

\bibitem{TTV}
D.K. Tomas, N. Tuneski and A. Vasudevarao, {\it Univalent Functions. A Primer}, De Gruyter, Berlin-Boston, 2018.

}

\end{thebibliography}
\end{document}